\documentclass[10 pt]{article}
\usepackage{amssymb,amsmath,amstext,amsgen,amsfonts,amsthm,fullpage}
\renewcommand{\epsilon}{\varepsilon}

\newcommand{\N}{\mathbb{N}}
\newcommand{\R}{\mathbb{R}}
\newcommand{\Z}{\mathbb{Z}}

\newcommand{\1}{{\bf 1}}
\newcommand{\supp}{\textbf{supp}}

\newtheorem{thm}{Theorem}[section]
\newtheorem{lem}{Lemma}[section]
\newtheorem{prop}{Proposition}[section]

\theoremstyle{remark}
\newtheorem{rmk}{Remark}[section]
\begin{document}

\title{Compression bounds for wreath products}
\author{Sean Li\footnote{Supported in part by NSF grants CCF-0635078 and CCF-0832795.}\\Courant Institute\\{\tt seanli@cims.nyu.edu}}
\date{}
\maketitle

\parindent = 0in

\begin{abstract}
We show that if $G$ and $H$ are finitely generated groups whose Hilbert compression exponent is positive, then so is the Hilbert compression exponent of the wreath $G \wr H$. We also prove an analogous result for coarse embeddings of wreath products. In the special case $G=\Z$, $H=\Z \wr \Z$ our result implies that the Hilbert compression exponent of $\Z\wr (\Z\wr\Z)$ is at least $1/4$, answering a question posed in \cite{ArGuSa,Te,NaPe2}.
\end{abstract}

\section{Introduction}

Given two groups, $G$ and $H$, we denote their wreath product $G \wr H$ to be the group of all pairs $(f,x)$ where $f : H \to G$ is finitely supported (i.e. $f(z) = e_G$ for all but finitely many $z \in H$) and $x \in H$.  We will let $\supp f$ denote $\{x \in H : f(x) \neq e_G\}$.  Let $T_xg(z) = g(x^{-1}z)$.  The product is given by
\begin{align*}
 (f,x) \cdot (g,y) = (z \mapsto f(z)g(x^{-1}z),xy) =: (fT_xg,xy)
\end{align*}
Using the usual heuristic of interpreting wreath products word metrics as lamplighter walks (see e.g. \cite{NaPe1}), we can interpret the metric as a traveling salesman tour on $H$.  For this paper, all groups will assumed to be infinite unless specified otherwise.  If $S$ and $T$ generate $G$ and $H$ respectively, then $G \wr H$ is generated by the set $\{(e_{G^H},t): t \in T\} \cup \{(\delta_s,e_H) : s \in S\}$ where $\delta_s$ is the function taking value $s$ at $e_H$ and $e_G$ everywhere else.  We will denote $\mathcal{L}_G(H)$ to be the wreath product where the generating set of $G$ is taken to be all elements. \\

Note that mappings $f$ of a metric space into $L_1$ induce pseudometrics given by
\begin{align*}
 d_f(x,y) = \|f(x) - f(y)\|_1.
\end{align*}
A cut in $G$ is a subset $A$ of $G$.  We say a cut $A$ separates a subset $B$ (denoted $A \vdash B$) if $B$ intersects both $A$ and $A^c$ nontrivially.  Given a measure $\mu$ on $\mathcal{C}$, a subset of the set of cuts of $G$, we can define a cut pseudometric by
\begin{align*}
 d_{\mu}(x,y) = \int_\mathcal{C} |\1_E(x) - \1_E(y)| ~d\mu(E).
\end{align*}
It is known that pseudometrics induced from embeddings into $L_1$ come from cut pseudometrics and vice versa \cite{DeLa,ChKl}. \\

For metric spaces $X$ and $Y$, we say that $f : X \to Y$ is a coarse embedding if there exist increasing unbounded functions $\rho_1,\rho_2 : \R_+ \to \R_+$ satisfying the inequality
\begin{align*}
 \rho_1(d_X(x,y)) \leq d_Y(f(x),f(y)) \leq \rho_2(d_X(x,y)).
\end{align*}

In \cite{GuKa}, Guentner and Kaminker gave a simple and elegant way of further strengthening coarse embeddability.  Let $G$ be a finitely generated group and $d$ the word metric with respect to its finite set of generators.  The $L_p$ compression exponent of $G$, which we denote by $\alpha_p^*(G)$, is the supremal $\alpha$ for which there exists a Lipschitz map $f : G \to L_p$ satsifying the inequality
\begin{align*}
 \|f(x) - f(y)\|_p \gtrsim d(x,y)^\alpha.
\end{align*}
We will use $\lesssim$ and $\gtrsim$ to denote inequality up to multiplication by some nonzero constant.  If $x \lesssim y$ and $y \lesssim x$, then we will say $x \asymp y$.  As all word metrics of a group with respect to finite generators are bi-Lipschitz equivalent, we see that the $L_p$ compression exponent of a group is in fact an algebraic invariant. \\

Furthermore, a map $f : G \to X$ where $X$ is a Banach space is called $G$-equivariant if it is given by the orbit of a vector $v \in X$ under the action of $G$ on $X$ by affine isometries.  That is, there is a representation $\alpha : G \to \text{Isom}(X)$ such that $f(g) = \alpha(g) \cdot 0$.  Note that the pseudometric $\|f(g) - f(h)\|_X$ is then $G$-invariant.  A group is said to have the Haagerup property if there exists an equivariant function $f : G \to L_2$ such that $\inf\{\|f(g) - f(h)\|_2 : d_G(g,h) \geq t\}$ tends to infinity with $t$.  We will say that this equivariant mapping is metrically proper. \\

We refer to the book \cite{ChCoJo} for more information on the Haagerup property, where in particular there is a discussion of which group operations preserve this property.  General semidirect products do not necessarily preserve the Haagerup property, or even the property of having positive compression exponent, as shown by the example $\Z^2 \rtimes \text{SL}_2(\Z)$.  Nevertheless, a recent breakthrough of de Cornulier, Stalder, and Valette \cite{CoStaVa} shows that wreath products do preserve the Haagerup property. \\

In this paper we show that the property of having a positive compression exponent is also preserved under wreath products.  Our approach uses crucially a tool from \cite{CoStaVa}, namely a method for lifting cut pseudometrics from a group $H$ to its associated lamplighter group $\mathcal{L}_G(H)$.  Several authors previously investigated the behavior of compression exponents under the wreath product operation in various special cases \cite{ArGuSa,AuNaPe,Ga,NaPe1,NaPe2,StVa,Te}.  All known examples with calculable compression exponent have shown that the wreath product operation preserves positivity of compression exponents, but it was unknown if this was true in general.  We prove that this is indeed the case when $p \in [1,2]$. \\

We first prove that positive $L_p$ compression exponent is preserved for the special case when $p = 1$. It is a standard fact that $L_p$ isometrically embeds into $L_1$ for $p \in (1, 2]$ (see e.g. \cite{BeLi}) so this reduction can always be done. We then lift the cut measure generated the cut decomposition of the $L_1$ embedding of $H$ to the cuts of $G \wr H$ using the method of \cite{CoStaVa}. By the correspondence between $L_1$ embeddings and cut measures, this gives us an embedding of $G \wr H$ into $L_1$, for which we show has positive $L_1$ compression.  To get the $L_p$ case, we then embed back into $L_p$ using a standard technique. This gives us the following theorem. 

\begin{thm}\label{thm1}
 For $p \geq 1$, we have: 
$$\alpha_p^*(G \wr H) \geq  \max\left\{\frac{1}{p},\frac{1}{2}\right\} \cdot \min\left\{\alpha_1^*(G), \frac{\alpha_1^*(H)}{1 + \alpha_1^*(H)}\right\}.$$
\end{thm}

It has been previously asked \cite{ArGuSa,Te,NaPe2} whether $\Z \wr (\Z \wr \Z)$ has a positive Hilbert compression exponent.  Using the fact that $\alpha_1^*(\Z \wr \Z) = 1$ \cite{NaPe1} we get a lower bound of $\alpha_2^*(\Z \wr (\Z \wr \Z)) \geq 1/4$.

\medskip

\begin{rmk}\label{rem:Q} Write $G_1=\Z\wr \Z$ and $G_{k+1}=\Z \wr G_k$. Since $G_1$ admits a bi-Lipschitz embedding into $L_1$~\cite{NaPe1}, it follows from Theorem~\ref{thm1} that for all integers $k$ we have $\alpha_1^*(G_k)\ge \frac{1}{k}$.  The question whether or not $\inf_{k\in \mathbb N} \alpha_1^*(G_k)$ is strictly positive remains an interesting open problem.
\end{rmk}

\medskip

In \cite{CoStaVa}, the authors showed that coarse embeddability into $L_p$ for $p \in [1,2]$ is preserved by wreath products.  We use the same construction as in the positive $L_p$ compression exponent theorem to give a somewhat simpler constructive proof with quantitative moduli bounds.

\begin{thm}\label{thm2}
 If $G$ and $H$ embed coarsely into $L_1$, then $G \wr H$ embeds coarsely into $L_p$ for $p \geq 1$.  Specifically, suppose we have mappings $\psi : H \to L_1$ and $\phi : G \to L_1$ satisfying compression bounds
 \begin{align*}
  &\rho_1(d(x,y)) \leq \|\psi(x) - \psi(y)\|_1 \leq \rho_2(d(x,y)), \\
  &\tau_1(d(x,y)) \leq \|\phi(x) - \phi(y)\|_1 \leq \tau_2(d(x,y)).
 \end{align*}
 Then there is a mapping $\Phi : G \wr H \to L_p$ satisfying the compression bounds
 \begin{align*}
  \nu_1(d((f,x),(g,y))) \lesssim \|\Psi(f,x) - \Psi(g,y)\|_p \lesssim \nu_2(d((f,x),(g,y))),
 \end{align*}
 where
 \begin{align*}
  \nu_2(t) &= t^{\max\left\{\frac{1}{2},\frac{1}{p}\right\}}, \\
  \nu_1(t) &= \left[\min\left\{ \sqrt{\min\left\{ \sum_{k=1}^n \rho_1(s_k) : s_k \geq 0, \sum_{k=1}^n s_k = t_1 \right\}} + \sum_{j=2}^m \tau_1(t_j) : t_j \geq 0, \sum_{j=1}^m t_j = t \right\}\right]^{\max\left\{\frac{1}{2},\frac{1}{p}\right\}}.
 \end{align*}
\end{thm}
In section 3, we show that $\nu_1$ is an increasing unbounded function. \\

We also show that equivariance is preserved by this construction.  Taken with coarse embeddability, we get a somewhat simpler proof of the Haagerup property result of \cite{CoStaVa},

\begin{thm}\label{thm3}
 If $G$ and $H$ have the Haagerup property, then so does $G \wr H$.  The compression bounds can be taken to be the same as those in Theorem 1.2.
\end{thm}

As a reduction to $L_1$ embeddings is part of the proof, we also show that the property of having a metrically proper equivariant maps into $L_1$ is preserved under wreath products.

\begin{rmk}
 As mentioned before, this method crucially uses the fact that $L_p$ isometrically embeds into $L_1$ for $p \in (1,2]$ to show that $L_p$ embeddability is preserved by wreath products. For $p > 2$, this is no longer true.  It is unknown if $G \wr H$ embeds coarsely into $L_p$ for $p > 2$ if $G$ and $H$ do so for $p > 2$ but not for $p \in [1,2]$.  The question of positive compression bounds also remains open, and it would be interesting to develop a method for this parameter range.
\end{rmk}

\subsection*{Acknowledgements} I am grateful to Assaf Naor for introducing me to the problems studied here and for suggesting that I look at~\cite{CoStaVa} in this context. I also thank him for  his improvement of my initial result, which was $\alpha_1^*(G \wr H)\ge  \min\left\{\alpha_1^*(G), \frac{\alpha_1^*(H)}{2}\right\}$. Note that, using the notation of Remark~\ref{rem:Q}, this weaker bound would yield the  bound $\alpha_1^*(G_k)\ge \frac{1}{2^{k-1}}$ rather than our result $\alpha_1^*(G_k)\ge \frac{1}{k}$, though, as asked in Remark~\ref{rem:Q}, it might be the case that there is actually a universal lower bound on $\alpha_1^*(G_k)$. Lastly, I am grateful to Assaf Naor for many constructive discussions and for his many revisions of this manuscript.

\section{Positive $L_p$ compression exponent}

We start by proving the following special case.

\begin{prop}
 $\alpha_1^*(\mathcal{L}_G(H)) \geq \frac{\alpha_1^*(H)}{1 + \alpha_1^*(H)}$.
\end{prop}

We begin by giving a method of extending measures on cuts of $H$ to measures on cuts of $\mathcal{L}_G(H)$. This is a specialization of the method used in \cite{CoStaVa}. Before we start the proof, as a warm-up, we give a description of the cuts for the case $C_2 \wr \Z^d$.  Note that the cuts in $\Z^d$ are just half-spaces orthogonal to some axis in $\Z^d$. The cut measure is just the discrete uniform measure on all the half-spaces. Elements of $C_2 \wr \Z^d$ can be thought of as finite subsets of $\Z^d$ with a special point for the initial lamplighter position. \\

Let $H$ be a half-space. Then, given a finite subset of points $A$ outside $H$, we define the cut $E(H, A)$ to be the collection of all elements of $C_2 \wr \Z^d$ whose initial lamplighter position is inside $H$ and whose subset outside of $H$ agrees with $A$. The total collection of cuts on $C_2 \wr \Z^d$ is simply the collection of all such cuts for all $H$ and all finite subsets outside of each $H$ . The cut measure is simply the discrete uniform measure on all the possible cuts. \\

We now give the full proof of Proposition 2.1. 

\begin{proof}
 As $\alpha_1^*(H) = \alpha$, for every $\epsilon \in (0,\alpha)$, we can find a Lipschitz function $\phi : H \to L_1$ so that for every $x,y \in H$, we have
 \begin{align*}
  \|\phi(x) - \phi(y)\|_1 \gtrsim d_H(x,y)^{\alpha - \epsilon}.
 \end{align*}
 
 Let $\rho$ be the cut measure on $\mathcal{C}$ induced by the cut decomposition of $\phi$.  We then know that $\rho$ induces a cut metric. We can complete $\mathcal{C}$ under taking complements and then deÞne a new cut measure $\mu(A) = \rho(A) + \rho(\check{A})$ where
 \begin{align*}
  \check{A} = \{ B^c : B \in A \}.
 \end{align*}
 The $\sigma$-algebra of the completion can simply be taken as the $\sigma$-algebra generated by $\mathcal{C}$ and $\check{\mathcal{C}}$ . We will still denote this $\mathcal{C}$ . Note then that $\mu$ is invariant under complements and differs from $\rho$ by at most a multiplicative factor of 2. Given a cut $B$ and finitely supported functions $f : B^c \to G$, we let 
 
 \begin{align*}
  E(B,f) &= \{(g,x) \in \mathcal{L}_G(H) : x \in B, ~g|_{B^c} = f\}.
 \end{align*}
 We will abuse notation and write $E(B,h)$ for $h : H \to G$ to mean the obvious cut using the restriction of $h$.  Note that even if there exists an element $y \in B$ such that $h(y) \neq h'(y)$, we still have that $E(B,h) = E(B,h')$ if $h|_{B^c} = h'|_{B^c}$. \\

 We can construct a measure $\nu$ on the set of cuts given by
 \begin{align*}
  \mathcal{H} = \bigcup_{B \in \mathcal{C}} \mathop{\bigcup_{f : B^c \to G}}_{\supp f < \infty} \{E(B,f)\}
 \end{align*}
 which induces a cut metric on $\mathcal{L}_G(H)$.  Indeed, consider the space
 \begin{align*}
  \tilde{\mathcal{C}} = \mathop{\bigsqcup_{f : H \to G}}_{\supp f < \infty} \{(B,f) : B \in \mathcal{C}\}.
 \end{align*}
 Note $\tilde{\mathcal{C}}$ is a countable union of sets bijective to $\mathcal{C}$.  Indeed, $H$ is countable and the set of finite subsets of $H$ is thus also countable.  Define the $\sigma$-algebra on $\tilde{\mathcal{C}}$ to be the $\sigma$-algebra generated by the $\sigma$-algebra of each of the factors.  There is an obvious measure, $\tilde{\mu}$ on $\tilde{\mathcal{C}}$ that restricts to $\mu$ on each of the factors.  That is, for every $f : H \to G$ and $A \subset \mathcal{C}$, we have that $\tilde{\mu}(\{(B,f) : B \in A\}) = \mu(A)$. \\

 Define the injection $\iota : \mathcal{H} \hookrightarrow \tilde{\mathcal{C}}$ given by $\iota(E(B, f )) = (B,f)$ where $f|_{B^c} = f$ and $f|_B = e_{G^B}$ . We can define a pullback $\sigma$-algebra on $\mathcal{H}$ so that $\iota$ is measurable. We then define the measure $\nu$ on $\mathcal{H}$ by $\nu(A) = \tilde{\mu}(\iota(A))$.  This is well defined as $\mathcal{H}$ has the pullback $\sigma$-algebra of $\iota$.  As $\iota$ is injective and $\mu$ is countably additive, $\nu$ is countably additive and thus indeed a measure.  \\

 One can see that
 \begin{align*}
  d_\nu((f,x),(g,y)) &= \nu(\{B \in \mathcal{H} : B \vdash \{(f,x), (g,y)\}\})\\ &\asymp \mu(\{B \in \mathcal{C} : B \vdash (\supp f^{-1}g) \cup \{x,y\}\}).
 \end{align*}
 This follows from easy case analyses.  In the $\gtrsim$ direction, suppose that $B \vdash (\supp f^{-1}g) \cup \{x,y\}$ (recall that we say $B$ separates $A$ or $B \vdash A$ if $A \cap B \neq \emptyset$ and $A \cap B^c \neq \emptyset$).  If $B \vdash \{x,y\}$ then from a case analysis of whether $x$ or $y$ is in $B$, either $E(B,f)$ or $E(B,g)$ separate $\{(f,x),(g,y)\}$.  No other possible cut can separate $\{(f,x),(g,y)\}$ using the same cut $B$.  Indeed if $E(B,h)$ is some other cut, then the restriction of $f$ and $g$ outside $B$ could not be equal to $h$ and so $(f,x)$ and $(g,y)$ would both not be in $E(B,h)$.  If $\{x,y\} \subset B$, then $E(B,f)$ and $E(B,g)$ both separate as there must be an element of $\supp f^{-1}g$ outside of $B$.  As before, no other $E(B,h)$ can separate.  If $\{x,y\} \cap B = \emptyset$, then $E(B^c,f)$ and $E(B^c,g)$ are the cuts that separate as $B$ contains an element of $\supp f^{-1}g$.  As $\mathcal{C}$ is closed under taking complements, we have that these cuts are in $\mathcal{H}$. \\

 In the $\lesssim$ direction, suppose $E(B,h)$ separates $\{(f,x),(g,y)\}$.  Note that $x$ and $y$ cannot both be in $B^c$ as then $(f,x)$ and $(g,y)$ are not in $E(B,h)$.  If $B \vdash \{x,y\}$, then we are done.  Suppose then that $x,y \in B$.  Then we must have that $g|_{B^c} \neq f|_{B^c}$ and so $(\supp f^{-1}g) \cap B^c \neq \emptyset$. \\

 In particular,
 \begin{align*}
  d_\nu((f,x),(g,y)) \gtrsim \max \left\{d_\mu(x,y) \cup \{d_\mu(x,z)\}_{z \in \supp f^{-1}g}\right\}.
 \end{align*}
 Thus, we have a function $F : \mathcal{L}_G(H) \to L_1$ such that
 \begin{align*}
  \|F(f,x) - F(g,y)\|_1 = d_\nu((f,x),(g,y)) \gtrsim \max \{d_\mu(x,y) \cup \{d_\mu(x,z)\}_{z \in \supp f^{-1}g}\}.
 \end{align*}

 This function is Lipschitz.  By the triangle inequality and the traveling salesman interpretation of the wreath product word metric, it suffices to check this for the cases when $(f,x)$ and $(g,y)$ differ by only a generator.  If $(f,x)$ and $(g,x)$ differ by an element of the form $(\delta_g,e_H)$ for $g \in G$, then $F(f,x) = F(g,x)$.  Indeed, $\supp f^{-1}g \cup \{x,x\} = \{x\}$, and there are no cuts that can separate a singleton.  Suppose then that $(f,x)$ and $(f,y)$ differ by $(0,t)$ where $t$ is a generator of $H$.  Then $\supp f^{-1}f \cup \{x,y\} = \{x,y\}$ and by above,
 \begin{align*}
  \|F(f,x) - F(f,y)\| \lesssim \mu(\{B \in \mathcal{C} : B \vdash \{x,y\}\}) = d_\mu(x,y) \lesssim d((f,x),(f,y)).
 \end{align*}
 By our choice of $\mu$, we have that
 \begin{align*}
  d((f,x),(g,y)) \gtrsim \|F(f,x) - F(g,y)\|_1 \gtrsim \max \{d_H(x,y)^{\alpha - \epsilon} \cup \{d_H(x,z)^{\alpha - \epsilon}\}_{z \in \supp f^{-1}g}\}.
 \end{align*}
 
 For a finitely supported function $f:H\to G$ define $\Lambda(f)\in \ell_1(H \times G)$ by 
 $$\Lambda(f)_{xy}=\left\{\begin{array}{ll}
 \frac12 & y=f(x),\\
 0 &\mathrm{otherwise.}
 \end{array}\right.$$
 Thus for $f,g:H\to G$ we have $\|\Lambda(f)-\Lambda(g)\|_1=|\supp f^{-1}g|$.  Now, consider the function
 \begin{align*}
  \Psi : \mathcal{L}_G(H) &\to \ell_1(H \times G)\oplus L_1 \\
  (f,x) &\mapsto \left(\frac{\alpha - \epsilon}{1 + \alpha - \epsilon} \Lambda(f)\right) \oplus \left(\frac{1}{1 + \alpha - \epsilon} F(f,x)\right).
 \end{align*}
 This function is clearly Lipschitz.  Bounding from below, we get
 \begin{align*}
  \|\Psi(f,x) - \Psi(g,y)\|_1 &\gtrsim \frac{\alpha - \epsilon}{1 + \alpha - \epsilon} |\supp f^{-1}g| + \frac{1}{1 + \alpha - \epsilon} \max \{d_H(x,y)^{\alpha - \epsilon} \cup \{d_H(x,z)^{\alpha - \epsilon}\}_{z \in \supp f^{-1}g}\} \\
  &\ge \max\left\{\left(|\supp f^{-1}g| \cdot \max \{d_H(x,y) \cup \{d_H(x,z)\}_{z \in \supp f^{-1}g}\}\right)^{\frac{\alpha - \epsilon}{1 + \alpha - \epsilon}}, \frac{d_H(x,y)^{\alpha - \epsilon}}{1+\alpha-\epsilon}\right\} \\
  &\gtrsim \left((1+|\supp f^{-1}g|) \cdot \max \{d_H(x,y) \cup \{d_H(x,z)\}_{z \in \supp f^{-1}g}\}\right)^{\frac{\alpha - \epsilon}{1 + \alpha - \epsilon}}\\
  \\
  &\gtrsim d((f,x),(g,y))^{\frac{\alpha - \epsilon}{1 + \alpha - \epsilon}}.
 \end{align*}
 The second inequality above is a consequence of the concavity of the logarithm, and the fourth inequality above comes from the triangle inequality.  Indeed, if $x = x_0,x_1,...,x_n=y$ is the shortest traveling salesman tour that starts from $x$, covers $\supp f^{-1}g$, and ends at $y$, then
 \begin{multline*}
  d((f,x),(g,y)) \lesssim \sum_{i=0}^{n-1} d_H(x_i,x_{i+1})\le \sum_{i=0}^{n-1} \left(d_H(x_i,x_0)+ d_H(x_{i+1},x_0)\right) =d_H(x,y)+2\sum_{z\in \supp f^{-1}g}d_H(x,z) 
   \\\lesssim (1+|\supp f^{-1}g|) \cdot \max\{d_H(x,y) \cup \{d_H(x,z)\}_{z \in \supp f^{-1}g}\}.
 \end{multline*}
 Taking $\epsilon$ to 0 finishes the proof.
\end{proof}

\begin{proof}[Proof of Theorem 1.1]
 Let $\alpha = \alpha_1^*(H)$ and $f : H \to L_1$ be a Lipschitz map with compression exponent $\alpha - \epsilon$ for some $\epsilon > 0$.  By the proposition, we then have a Lipschitz map $\Psi : \mathcal{L}_G(H) \to L_1$ with compression exponent $\frac{\alpha - \epsilon}{1 + \alpha - \epsilon}$.  Thus, we have proven that $\alpha_1^*(\mathcal{L}_G(H)) \geq \frac{\alpha_1^*(H)}{1 + \alpha^*(H)}$.  Note that as $\alpha_1^*(\mathcal{L}_G(H)) \leq 1$ always, Theorem 3.3 of \cite{NaPe1} gives that $\alpha_1^*(G \wr H) \geq \min\left\{\alpha_1^*(G),\frac{\alpha_1^*(H)}{1 + \alpha_1^*(H)}\right\}$. \\

The case  $p> 1$ in Theorem~\ref{thm1} is  a consequence of the general fact that $\alpha_p^*(\Gamma)\ge \max\left\{\frac{1}{p},\frac{1}{2}\right\}\alpha^*_1(\Gamma)$, which holds for any finitely generated group $\Gamma$. This simple fact is explained in~\cite{NaPe2} (see the paragraph following question 10.3 there for a more general statement). For completeness we will now recall how this is proved.  As there is a map $T_p : L_1 \to L_p$ such that $\|T_p(x) - T_p(y)\|_p^p = \|x - y\|_1$ \cite{Ma,WeWi} we can compose the embedding of $\Gamma$ into $L_1$ with $T_p$ to get the result that $\alpha_p^*(\Gamma) \geq \frac{1}{p} \alpha_1^*(\Gamma)$.  For $p > 2$, we can first embed into $L_2$ and then use the fact that $L_2$ embeds isometrically into $L_p$ for all $p$ (see e.g. \cite{Wo,BeLi}) to get the final bound $\alpha_p^*(\Gamma) \geq \max\left\{\frac{1}{2},\frac{1}{p}\right\} \alpha_1^*(\Gamma)$.  Applying this to our lower bound for $\alpha_1^*(G \wr H)$ gives the proof of Theorem 1.1.
\end{proof}

\section{Coarse embeddability}

 We will prove only the $L_1$ embedding case as the general $L_p$ case follows from the second half of the proof of Theorem 1.1.  Note that as word metrics take only integer values, we may view the bounding functions as unbounded increasing functions from $\Z_+$ to $\R_+$.  We can assume 0 is mapped to 0. \\

We begin by proving a quantitative analogue of a lemma from \cite{NaPe1}.  The proof will be mostly identical to the original version with a few key changes.

\begin{lem}
 If $G$ and $\mathcal{L}_G(H)$ both coarsely embed into $L_1$ then $G \wr H$ coarsely embeds into $L_1$.
\end{lem}

\begin{proof}
 Let $\ell_1(H,G;\text{fin})$ denote the metric space of all finitely supported functions $f : H \to G$ equipped with the metric
 \begin{align*}
  d_{\ell_1}(f,g) = \sum_{z \in H} d_G(f(z),g(z)).
 \end{align*}

 Then one sees that
 \begin{align*}
  d_{G \wr H}((f,x),(g,y)) \asymp d_{\mathcal{L}_G(H)}((f,x),(g,y)) + d_{\ell_1}(f,g).
 \end{align*}

 Indeed, we may suppose $(g,y) = ({\bf e},e)$ as the metrics are $G \wr H$-invariant.  Then to move from $({\bf e},e)$ to $(f,x)$ is the same as visiting all locations of $\supp f$ and at each location moving $G$ from $e$ to $f(z) \in G$. \\

 Let $\phi : G \to L_1$ and $\Psi : \mathcal{L}_G(H) \to L_1$ be coarse embeddings with the bounds
 \begin{align*}
  \tau_1(d(x,y)) &\leq \|\phi(x) - \phi(y)\|_1 \leq \tau_2(d(x,y)), \\
  \xi_1(d_{\mathcal{L}_G(H)}((f,x),(g,y))) &\leq \|\Psi(f,x) - \Psi(g,y)\|_1 \leq \xi_2(d_{\mathcal{L}_G(H)}((f,x),(g,y))).
 \end{align*}

 Define the function $F : G \wr H \to L_1 \oplus \ell_1(H,L_1;\text{fin})$ by
 \begin{align*}
  F(f,x) = \Psi(f,x) \oplus (\phi \circ f).
 \end{align*}

 We have that $d_{G \wr H}((f,x),(g,y)) \asymp d_{\mathcal{L}_G(H)}((f,x),(g,y)) + d_{\ell_1(H,G)}(f,g)$.
 \begin{align*}
  \|F(f,x) - F(g,y)\|_1 &= \|\Psi(f,x) - \Psi(g,y)\|_1 + \sum_{z \in H} \|\phi(f(z)) - \phi(g(z))\|_1 \\
  &\leq \xi_2(d_{\mathcal{L}_G(H)}((f,x),(g,y))) + \sum_{z \in H} \tau_2(d(f(z),g(z))) \\
  &\leq \eta_2\left(d_{\mathcal{L}_G(H)}((f,x),(g,y)) + \sum_{z \in H} d(f(z),g(z))\right) \\
  &\lesssim \eta_2(d_{G \wr H}((f,x),(g,y)))
 \end{align*}
 where $\eta_2(t) := t \cdot (\xi_2(t) + \tau_2(t))$.  As $\eta_2$ is clearly unbounded increasing, we have an upper bound.  Bounding from below, we have
 \begin{align*}
  \|F(f,x) - F(g,y)\|_1 &\geq \xi_1(d_{\mathcal{L}_G(H)}((f,x),(g,y))) + \sum_{z \in H} \tau_1(d(f(z),g(z))) \\
  &\geq \eta_1\left(d_{\mathcal{L}_G(H)}((f,x),(g,y)) + \sum_{z \in H} d(f(z),g(z))\right) \\
  &\geq \eta_1(d_{G \wr H}((f,x),(g,y))).
 \end{align*}
 where
 \begin{align*}
  \eta_1(t) := \min \left\{\xi_1(t_1) + \sum_{i=2}^n \tau_1(t_i) : t_i \geq 0, \sum_{i=1}^n t_i = t\right\}
 \end{align*}

 One can think of this function as evaluated on integer partitions of $t$.  To show that $\eta_1$ is increasing, take a partition of $t+1$.  If all the elements in the partition are of size 1, then the value would be $\xi_1(1) + \sum_{i=1}^t \tau_1(1)$ which is greater than the value for the partition of $t$ of all sizes 1 as $\xi_1(1) > 0$ and $\tau_1(1) > 0$.  If the partition $t + 1$ has an element of value greater than 1, then reducing this by one gives a partition of $t$ of lesser value as $\xi_1$ and $\tau_1$ are increasing.  Thus, we see that the minimum of the partitions of $t+1$ always is greater than the minimum of the partitions of $t$.  It remains to show that $\eta_1$ is unbounded. \\

 Let $M > 0$.  By rescaling, we may suppose that $\xi_1(1) \geq 1$ and $\tau_1(1) \geq 1$.  Let $N > 0$ such that $\tau_1(N) \geq M$ and $\xi_1(N) \geq M$.  Consider the possible partitions of $MN$.  If there are more than $M$ elements of the partition, then as $\tau_1(1) \geq 1$ and $\xi_1(1) \geq 1$, the summation associated to this partition would have value greater than $M$.  Thus, the number of elements of the partition has to be less than $M$.  However, one of elements in the partitions has to have value greater than $N$ by pigeonhole principle and so either of the $\xi_1$ or $\tau_1$ value of this element is greater than $M$.  Thus, $\eta_1(MN) > M$ and so $\eta_1$ is unbounded. \\

 Note that as $\xi_2$ and $\tau_2$ are increasing unbounded functions, $\eta_2$ grows superlinearly.  However, as $F$ is Lipschitz, we can always use the triangular inequality to give a linear upper bound based on the expansion between generators.  Thus, we may replace $t \cdot (\zeta_2(t) + \tau_2(t))$ with a linear compression bound.
\end{proof}

\begin{proof}[Proof of Theorem 1.2]
 By the preceding remarks and lemma, it suffices to show that $\mathcal{L}_G(H)$ embeds coarsely into $L_1$ when $G$ and $H$ embed coarsely into $L_1$.  As such, let $\psi : H \to L_1$ be an embedding with bounds
 \begin{align*}
  \rho_1(d(x,y)) \leq \|\psi(x) - \psi(y)\|_1 \leq \rho_2(d(x,y)).
 \end{align*}

 Using the cut decomposition, we construct the function $F : \mathcal{L}_G(H) \to L_1$ from the measure $\nu$ on cuts of $\mathcal{L}_G(H)$ as before.  As before, we have that
 \begin{align*}
  \|F(f,x) - F(g,y)\|_1 \asymp \mu(\{B \in \mathcal{C} : B \vdash (\supp f^{-1}g) \cup \{x,y\}\}).
 \end{align*}
 In particular,
 \begin{align*}
  \|F(f,x) - F(g,y)\|_1 \gtrsim \max \{\rho_1(d(x,y)) \cup \{\rho_1(d(x,z))_{z \in (\supp f^{-1}g) \cup \{y\}}\}\},
 \end{align*}
 and
 \begin{align*}
  \|F(f,x) - F(g,y)\|_1 &\lesssim \mu(\{B \in \mathcal{C} : B \vdash (\supp f^{-1}g) \cup \{x,y\}\}) \\
  &\leq \sum_{u,v \in (\supp f^{-1}g) \cup \{x,y\}} \|\psi(u) - \psi(v)\|_1 \\
  &\leq \sum_{u,v} \rho_2(d(u,v))
 \end{align*}

 Construct the mapping $\Psi(f,x) := \Lambda(f) \oplus F(f,x)$.  Bounding from above, we have that
 \begin{align*}
  \|\Psi(f,x) - \Psi(g,y)\|_1 &\lesssim |\supp f^{-1}g| + \sum_{u,v \in (\supp f^{-1}g) \cup \{x,y\}} \rho_2(d(u,v)) \\
  &\leq d((f,x),(g,y)) + \sum_{u,v} \rho_2(d((f,x),(g,y))) \\
  &\leq \tau_2(d((f,x),(g,y))).
 \end{align*}
 where $\tau_2(t) := t + t^2 \cdot \rho_2(t)$.  Bounding from below, we have
 \begin{align*}
  \|\Psi(f,x) - \Psi(g,y)\|_1 &\gtrsim |\supp f^{-1}g| + \max \{\rho_1(d(x,y)) \cup \{\rho_1(d(x,z))_{z \in (\supp f^{-1}g) \cup \{y\}}\}\} \\
  &\gtrsim \max\left\{\sqrt{|\supp f^{-1}g| \cdot \max \{\rho_1(d(x,y)) \cup \{\rho_1(d(x,z))_{z \in (\supp f^{-1}g) \cup \{y\}}\}\}}, \rho_1(d(x,y))\right\} \\
  &\gtrsim \sqrt{(1 + |\supp f^{-1}g|) \cdot \max \{\rho_1(d(x,y)) \cup \{\rho_1(d(x,z))_{z \in (\supp f^{-1}g) \cup \{y\}}\}\}} \\
  &\geq \tau_1(d((f,x),(g,y))).
 \end{align*}
 where
 \begin{align*}
  \tau_1(t) := \sqrt{\min \left\{\sum_{k=1}^n \rho_1(s_k) : s_k \geq 0, \sum_{k=1}^n s_k = s\right\}}.
 \end{align*}

 The second inequality came from using the AM-GM inequality.  We can also use Young's inequality rather than the AM-GM inequality to improve the lower bound for certain lower moduli.  Using a similar proof as above, we can see that $\tau_1$ is an increasing unbounded function.  Thus, we have that $\mathcal{L}_G(H)$ embeds coarsely into $L_1$.  Composing these compression bounds with those of Lemma 2.1 gives us the necessary bounds.
\end{proof}

\section{The Haagerup property}

In this section, we use the ideas from above to show that equivariant maps into Hilbert space can be amalgamated to give an equivariant map defined on $G \wr H$.  As before we first prove the $L_1$ analogue.

\begin{thm}
 If $G$ and $H$ admit metrically proper equivariant mappings into $L_1$, then so does $G \wr H$.
\end{thm}

As above, we need to prove the following lemma

\begin{lem}
 If $G$ and $\mathcal{L}_G(H)$ admit equivariant mappings into $L_1$, then so does $G \wr H$.
\end{lem}

\begin{proof}
 Let $\psi : G \to L_1$ and $\phi : \mathcal{L}_G(H) \to L_1$ be equivariant maps with associated actions $\tau$ and $\pi$.  We would like to show that
 \begin{align*}
  \Psi : G \wr H &\to L_1 \oplus \ell_1(H,L_1;\text{fin}) \\
  (f,x) &\mapsto \phi(f,x) \oplus (\psi \circ f)
 \end{align*}
 is equivariant.  We will express elements of $\ell_1(H,X;\text{fin})$ as elements of the direct product $\bigoplus_{h \in H} X$.  The action of $H$ on this direct product is then permutation of coordinates.  This is precisely the action of $T_h$ for $h \in H$.  Note that the semidirect product of $\bigoplus_{h \in H} G_h$ with $H$ by this action is just the wreath product $G \wr H$.  Consider the group action of $G \wr H$ on $L_1 \oplus \ell_1(H,L_1;\text{fin})$
 \begin{align*}
  \theta(f,x) \left(u,\sum_{h \in H} g_h\right) = \left( \pi(f,x)u, \sum_{h \in H} \tau(f(h)) g_{x^{-1}h} \right)
 \end{align*}
 where $g_h \in L_1$.  It is then straightforward from the formulas to see that
 \begin{align*}
  \Psi((f,x) \cdot (g,y)) = \theta(f,x) \cdot \Psi(g,y) + \Psi(f,x).
 \end{align*}
\end{proof}

We also require the following theorem \cite{La} (see also the exposition in \cite{FlJa}).

\begin{thm}[Lamperti's Theorem]
 Let $U$ be an isometry of $L_1$ onto itself.  Then there is a Borel measurable self-mapping $\varphi$ of $[0,1]$ that is bijective almost everywhere and a $u \in L_1$ such that
 \begin{align*}
  U\psi = u \cdot (\psi \circ \varphi).
 \end{align*}
 Furthermore,
 \begin{align*}
  \int_{\varphi^{-1}(E)} |u| ~dt = \int_E ~dt.
 \end{align*}
 for every Borel set $E$.
\end{thm}

As $U$ is an isometry, it is clear that $\varphi$ cannot map a set of positive measure to a set of measure 0 and vice versa.  In addition, $u$ must be nonzero on a set of full measure.

\begin{proof}[Proof of Theorem 4.1]
 By the preceding lemma, it suffices to show that the map $F$ on $\mathcal{L}_G(H)$ constructed as before from $\psi : H \to L_1$ is equivariant.  Indeed, the mapping $(f,x) \mapsto \Lambda(f)$ is equivariant and so the entire embedding $\Psi(f,x) = F \oplus \Lambda(f)$ would be equivariant.  As we are using the same construction of $\Psi$ as before, that $\Psi$ is metrically proper will follow from the arguments of the Section 3. \\

 Recall that the cuts generated in the decomposition of $\psi$ are given by the cut map
 \begin{align*}
  S : [0,1] \times \R &\to \mathcal{C} \\
  (y,t) &\mapsto \{h \in H : t^{-1} \cdot \psi(h)(y) > 1\},
 \end{align*}
 and the measure $\rho$ is the pushforward of the Lebesgue measure on $[0,1] \times \R$ by $S$.  As before, we complete $\mathcal{C}$ under taking complements and define a new complement invariant measure $\mu$ on $\mathcal{C}$ from $\rho$.  Having defined $\mu$, we extend it to $\nu$ on $\mathcal{H}$ as before and from this we get a map $F : \mathcal{L}_G(H) \to L_1(\mathcal{H},\nu)$ with the desired properties. \\

 We would like to show that there exists an isometric group action of $\mathcal{L}_G(H)$ on $L_1(\mathcal{H},\nu)$ such that $F(gh) = \pi(g) \cdot F(h) + F(g)$.  We accomplish this by showing that the natural action of $H$ on the set of cuts $\mathcal{C}$ by left multiplication is measure preserving.  Note that we can extend the action of $H$ to $\mathcal{H}$ by
 \begin{align*}
  hE_i(B,f) = E_i(hB, T_hf).
 \end{align*}
 The problem comes from determining whether $hB$ is in $\mathcal{C}$.  Given a finitely supported function $g : H \to G$, we can also specify the action
 \begin{align*}
  gE_i(B,f) = E_i(B,g|_{B^c} \cdot f).
 \end{align*}
 These two actions are easily seen to be compatible with the group operation.  Thus, to show that this is actually a group action, it suffices to show that $\mathcal{C}$ is $H$-invariant except possibly on a set of measure 0.  To show that this group action is isometric, we require that $\nu$ be $G \wr H$-invariant.  By the pullback construction of $\nu$, it suffices to show that $\rho$ is $H$-invariant. \\

 Let $\varphi$ and $u$ be the functions associated to $\pi(h)$ by Lamperti's theorem.  From equivariance, we have that
 \begin{align*}
  \psi(hg)(\varphi^{-1}(y)) = \pi(h) \cdot \psi(g)(\varphi^{-1}(y)) + \psi(h)(\varphi^{-1}(y)) = u(\varphi^{-1}(y)) \cdot \psi(g)(y) + \psi(h)(\varphi^{-1}(y))
 \end{align*}
 for almost every $y \in [0,1]$.  As we only need to prove $H$-invariance on a full measure subset of $\mathcal{C}$, we may suppose that $x = \varphi^{-1}(y)$ is defined and $a = u(x)$ and $b = \psi(h)(x)$ are finite.  Suppose $t > 0$.  If $a \cdot (at + b) > 0$, we have that $hS(y,t) = S(x, at + b)$.  If $a \cdot (at + b) < 0$ then $hS(y,t) = S(x, at + b)^c$.  As we required that $\mathcal{C}$ be closed under complement, this is not a problem.  The case when $t < 0$ can also be similarly analyzed.  As $\varphi$ cannot take measure 0 sets onto sets of positive measure and vice versa, this shows that the set of cuts is $H$-invariant up to a set of cuts of measure 0.  It remains to show that $\rho$ is $H$-invariant. \\

 As before, fix $h \in H$ and let $\varphi$, $u$ be the functions associated to the isometry $\pi(h)$.  Let $B \in \mathcal{C}$.  Note by the cut map that the quantity $\mu(\{B\}) = \rho(\{B\}) + \rho(\{B^c\})$ can be thought of as the area bound between the two family of curves $\{\psi(g) : g \in B\}$ and $\{\psi(g) : g \in B^c\}$, that is the Lebesgue measure of the points of $[0,1] \times \R$ that is below all the graphs of one family and above all the graphs of the other. \\
 
 Note that $\pi(h)$ induces a self-mapping of $[0,1] \times \R$.  Indeed,  taking the arguments of the $H$-invariance into account, the transformation can be given by
 \begin{align*}
  \pi(h)(x,y) = (\varphi^{-1}(x), y \cdot u(\varphi^{-1}(x)) + \psi(h)(\varphi^{-1}(x))).
 \end{align*}

 By Lamperti's theorem, this transformation, which is defined on a set of full measure, is precisely the one that takes the graph of $\psi(g)$ to the graph of $\psi(hg)$.  Note that vertical ordering of the graphs is preserved by this transformation at each $x \in [0,1]$ except with the possibility of a flip.  Let $E \times F \subset [0,1] \times \R$ be a measurable subset.
 \begin{align*}
  \int_{\pi(h)(E \times F)} ~dx ~dy = \int_{\varphi^{-1}(E)} \int_F |u(x)| ~dy ~dx = |F| \int_{\varphi^{-1}(E)} |u| ~dx = |E||F|
 \end{align*}

 Thus, we see that $\pi(h)$ is measure preserving for the set of generators of the $\sigma$-algebra of $[0,1] \times \R$.  As area is preserved, we have that $\mu(hB) = \mu(B)$.  It follows that $H$ induces a measure preserving transformation on $(\mathcal{H},\nu)$ and subsequently, an isometry of $L_1(\mathcal{H},\nu)$, which we will still denote $\pi$. \\

 Note that the cut map $S : [0,1] \times \R \to (\mathcal{C},\mu)$ induces an isometric embedding of $L_1(\mathcal{C},\mu)$ into $L_1([0,1] \times \R)$.  It is readily seen that $L_1(\tilde{C},\tilde{\mu})$ is a countable $\ell_1$ sum of $L_1(\mathcal{C},\mu)$ and so isometrically embeds into $\left(\sum_{j=1}^\infty L_1([0,1] \times \R)\right)_1$ which is isometric to $L_1$.  Note that $\iota$, the injection of $\mathcal{H}$ into $\tilde{\mathcal{C}}$, has the property that $\nu(\iota^{-1}(E)) \leq \tilde{\mu}(E)$.  This gives that the induced map $\iota^* : L_1(\tilde{\mathcal{C}},\tilde{\mu}) \to L_1(\mathcal{H},\nu)$ is continuous and onto and so $L_1(\mathcal{H},\nu)$ is separable. \\

 As $L_1(\mathcal{H},\nu)$ is separable, we know that it is isometric to one of the following spaces \cite{Wo}:
 \begin{align*}
  L_1, ~~~\ell_1, ~~~\{\ell_1^n\}_{n=1}^\infty, ~~~L_1 \oplus \ell_1, ~~~\{L_1 \oplus \ell_1^n\}_{n=1}^\infty.
 \end{align*}

 If $L_1(\mathcal{H},\nu)$ is isometric to $L_1$, then we are done.  Otherwise, we need to embed $L_1(\mathcal{H},\nu)$ into $L_1$ and define a suitable isometric action.  The argument will follow closely to the ones made in \cite{NaPe2}. \\
 
 If $L_1(\mathcal{H},\mu)$ is isometric to $\ell_1$, then Lamperti's theorem tells us that $\pi(g)e_i = \theta_i^g e_{\tau^g(i)}$ for all $i \in \N$ where $\{e_i\}_{i=1}^\infty$ are the standard coordinate basis for $\ell_1$, the function $\tau^g: \N \to \N$ is bijective, and $|\theta^g_i| \equiv 1$.  Embedding $L_1(\mathcal{H},\nu)$ into $L_1$ by the standard mapping $\varphi : x \mapsto \sum_{i=1}^\infty 2^i x_i \chi_{[2^{-i},2^{-i+1}]}$, we can define the isometric action of $H$ on $L_1$, which we will still denote $\pi$, by
 \begin{align*}
  \pi(g)f(t) := \theta_i^g f(2^{i-\tau^g(i)} t).
 \end{align*}
 It is immediate to check that $\pi$ and $\varphi \circ f$ satisfy the necessary equivariance relation. \\

 For the cases when $L_1(\mathcal{H},\nu)$ is isometric to $L_1 \oplus \ell_1(S)$ where $S$ is a countable set, we use Lamperti's theorem to show that isometric automorphisms map disjoint functions to disjoint functions and indicators of atoms to indicators of atoms.  Thus, $\pi(g)$ restricts to isometries of each summand.  By embedding $L_1 \oplus \ell_1(S)$ into $L_1 \oplus L_1$ (which is isometric to $L_1$), we can get the necessary isometric action by using the results from above on each of the $L_1$ summands.
\end{proof}

\begin{proof}[Proof of Theorem 1.3]
 We will use the same embedding as above except for a change in the first step of embedding $L_2$ into $L_1$.  Given an equivariant map $\psi : H \to L_2$, the proof of lemma 2.3 in \cite{NaPe1} shows that there exists an equivariant map $\tilde{\psi} : H \to L_1$ such that $\|\tilde{\psi}(x) - \tilde{\psi}(y)\|_1 = C \cdot \|\psi(x) - \psi(y)\|_2$ for some constant $C$. \\

 In conjunction with the previous theorem, we have that the existence of an equivariant map $f : H \to L_2$ gives an equivariant map $F : G \wr H \to L_1$.  Let $\pi$ be the associated group action on $L_1$.  It remains to show that the map $T \circ F : G \wr H \to L_2$ is also equivariant.  Recall that the embedding $T$ of $L_1[0,1]$ into $L_2([0,1] \times \R)$ can be thought of mapping $f$ to $(x,y) \mapsto 1 - \1_{[0,y]}(f(x))$.  This maps a function on $[0,1]$ to the characteristic function on the area bound by its graph.  As above, we may think of $\pi(f,x)$ as a measure preserving transformation of $[0,1] \times \R$ given by
 \begin{align*}
  \pi(f,x)(y,t) = (\varphi^{-1}(y), t \cdot u(\varphi^{-1}(y)) + F(f,x)(\varphi^{-1}(y))).
 \end{align*}

 This induces an isometry of $L_2([0,1] \times \R)$.  Given $\pi(f,x)$ and $\pi(g,y)$ with corresponding functions, $u_1,\varphi_1$ and $u_2,\varphi_2$, we see that
 \begin{align*}
  \pi((f,x) (g,y)) \cdot \psi(t) = \pi(f,x) (u_2(t) \cdot \psi(\varphi_2(t))) = u_1(t) \cdot u_2(\varphi_1(t)) \cdot \psi(\varphi_2 \circ \varphi_1(t)).
 \end{align*}
 Thus, the associated functions for $\pi((f,x) \cdot (g,y))$ are $u_1(t) \cdot u_2(\varphi_1(t))$ and $\varphi_2 \circ \varphi_1$.  The first coordinate of the mapping clearly preserves the group structure.  For the second coordinate, using the fact that $F$ is equivariant with respect to the action $\pi$, we get
 \begin{align*}
  \pi(f,x)\pi(g,y) \cdot (y,t) &= \pi(f,x)(\varphi_2^{-1}(y), t \cdot u_2(\varphi_2^{-1}(y)) + F(g,y)(\varphi_2^{-1}(y))) \\
  &= (\varphi_1^{-1}(\varphi_2^{-1}(y)), t \cdot u_2(\varphi_2^{-1}(y)) u_1(\varphi_1^{-1}(\varphi_2^{-1}(y))) \\
  &~~~+ u_1(\varphi_1^{-1}(\varphi_2^{-1}(y)))F(g,y)(\varphi_2^{-1}(y)) + F(f,x)(y')) \\
  &= (y', t \cdot u_2(\varphi_1(y')) u_1(y') + F((f,x)(g,y))(y')) \\
  &= \pi((f,x)(g,y))(y,t).
 \end{align*}
 Thus, $\pi$ is indeed a group action.  By the same arguments as above, this is a measure preserving group action and so induces a group action of $L_2([0,1] \times \R)$ which we denote $\tau$.  One can see from the definition of $T$ then that $T \circ F((f,x) \cdot (g,y)) = \tau(f,x) \cdot (T \circ F(g,y)) + T \circ F(f,x)$.
\end{proof}

\bibliographystyle{amsalpha}
\bibliography{l1wreath}
\end{document}